\documentclass[a4paper,10pt]{amsart}

\usepackage{amsmath}
\usepackage{amsfonts}
\usepackage{amssymb}
\usepackage{graphicx}
\usepackage{color}
\usepackage{soul,xcolor} 
\usepackage{cite}

\usepackage{amsthm}
\usepackage[all,cmtip]{xy}

\newcommand{\mbR}{\mathbb{R}}
\newcommand{\mbC}{\mathbb{C}}

\newcommand{\mbP}{\mathbb{P}}

\DeclareMathOperator{\Hom}{Hom}

\DeclareMathOperator{\Stab}{Stab}

\DeclareMathOperator{\coh}{coh}

\newcommand*{\coloneq}{\mathrel{\mathop:}=}

\theoremstyle{plain}
\newtheorem{theorem}{Theorem}[section]

\newtheorem{lemma}[theorem]{Lemma}

\theoremstyle{definition}
\newtheorem{definition}[theorem]{Definition}

\theoremstyle{remark}
\newtheorem{remark}[theorem]{Remark}

\usepackage{hyperref}

\hypersetup{
	colorlinks=true,         
	linkcolor=blue,          
	citecolor=blue,          
	hyperfootnotes=false,
}

\begin{document}
\bibliographystyle{alpha}

\title[Why care about Bridgeland stability conditions?]{Why should a birational geometer care about Bridgeland stability conditions?}

\author{Claudio Fontanari}
\address{Claudio Fontanari \\ Universit\`a degli Studi di Trento, Dipartimento di Matematica, via Sommarive 14, 38123 Povo, Trento, Italy.}

\email{claudio.fontanari@unitn.it} 

\author{Diletta Martinelli}
\address{Diletta Martinelli \\ Department of Mathematics, Imperial College London, 180 Queen's Gate, London SW7 2AZ, UK.}
\email{d.martinelli12@imperial.ac.uk}

\thanks{This research was partially supported by FIRB 2012 "Moduli spaces and Applications" and by GNSAGA of INdAM (Italy).}

\begin{abstract} 
In this survey we borrow from Coskun and Huizenga an example of application of Bridgeland stability conditions 
to birational geometry and we rephrase it without assuming any previous knowledge about derived categories.
\end{abstract}

\date{\today}

\maketitle

\section{Introduction}

\noindent
Before addressing the question in the title, perhaps we first need to justify why on earth the answer 
should come from two newcomers into the Bridgeland world. Indeed, while approaching the somehow 
exotic land of derived categories, we deeply felt the urgency of a strong motivation rooted in classical     
algebraic geometry. Even if it is well known that applications of Bridgeland stability conditions to birational
geometry are many and fruitful, it might seem (so it was for us) that in order to appreciate 
the geometric content of the theory one needs to be already involved in its jungle of technicalities.  
However, in the end we realized that at least some geometric example could be presented with only 
a little amount of machinery. Now we would like to humbly share the results of our efforts to get the 
point, without any claim of originality and exhaustiveness. 

\noindent
Let $X$ be a connected projective scheme over an algebraically closed field $k$ of characteristic zero. 
If we fix an ample line bundle $\mathcal{O}_X(1)$ on $X$ and a polynomial $P \in \mathbb{Q}[z]$, 
then according to \cite{HL}, Theorem 4.3.4, there is a projective scheme whose closed points 
are in bijection with $S$-equivalence classes of Gieseker semistable sheaves with Hilbert polynomial $P$.
In particular, we can consider the moduli space $M = M_{\mathbb{P}^2}(r, c_1,c_2)$ of $S$-equivalence 
classes of semistable sheaves of rank $r$ and Chern classes $c_1$ and $c_2$ on the projective plane 
(see for instance \cite{LP}, Part II). 

\noindent
As explained in the preface to \cite{HL}, there are good reasons to study moduli spaces of sheaves. 
In particular, \emph{they provide examples of higher dimensional algebraic varieties with a rich and 
interesting geometry. In fact, in some regions in the classification of higher dimensional varieties the
only known examples are moduli spaces of sheaves on a surface.} From a birational geometry perspective, 
many natural questions arise: describe the ample cone, determine the effective cone, run an MMP, 
and so on. 

\noindent
Quite recently, remarkable progress in the field has been obtained by the application of 
Bridgeland stability conditions introduced in \cite{B1}. In particular, Bayer and 
Macr\`i have described the nef cone of the moduli space of Gieseker stable sheaves on a K3 surface 
in \cite{BM2} and Coskun and Huizenga have computed the ample cone of Gieseker semistable 
sheaves on $\mathbb{P}^2$ in \cite{CH}.

\noindent
Here we have chosen the description of the ample cone as our guiding example in order to motivate 
an ideal reader with a background in classical algebraic geometry and without any previous knowledge 
about derived categories. In this spirit, we are going to adopt a slightly unconventional order: indeed, 
we start with geometric applications, by taking stability conditions as a sort of black box, and we turn 
to algebraic formalism only at the end. More precisely, in Section \ref{sta} we present the beautiful 
geometry related to Bridgeland stability conditions. Next, in Section \ref{amp} we focus on the case 
of $\mathbb{P}^2$ and we outline the procedure to compute the ample cone of $M$ introduced in 
\cite{CH}. Finally, in Section \ref{def} we collect precise definitions in the case of $\mathbb{P}^2$ 
and specific references for the general case. 

\subsubsection*{Acknowledgements}
\noindent
This project started during the visit of the second--named author at the Department of Mathematics 
at the University of Trento. She wishes to thank this institute for the warm hospitality. 
Both authors would like to thank Arend Bayer for useful comments and suggestions.

\section{Geometry of stability conditions}\label{sta}
\noindent
Let $\mathcal{A} = \coh(X)$ be the abelian category of coherent sheaves on a smooth projective variety $X$ defined over the complex numbers.
Bridgeland's key idea was to introduce stability conditions not on the abelian category of coherent sheaves, but on the bounded derived category $D^b(\mathcal{A})$. 
As anticipated in the introduction, for now we consider the notion of stability condition as a black box. However, in order to avoid cheating too much, it may be useful 
to have in mind at least the definition of derived category, which we now recall from \cite{B}, Section 3.

\noindent
Let $C^b(\mathcal{A})$ be the category of bounded complexes: objects are complexes 
\begin{equation*}
E^{\bullet} = \ldots \to E^i \stackrel{d^i}{\to}  E^{i+1} \stackrel{d^{i+1}}{\to} E ^{i+2} \to \ldots
\end{equation*}
with $d^{i+1}\circ d^i = 0$ and $H^i(E) = 0$ for all but finitely many $i$, and morphisms $f^{\bullet}: E^{\bullet} \to F^{\bullet}$ are morphisms $f_i: E^i \to F^i$ that commute with the differential.
A morphism $f^{\bullet}$ is called a quasi-isomorphism if the induced morphism in cohomology $f_* : H^i(E^{\bullet}) \to H^i(F^{\bullet})$ is an isomorphism for all $i$.

\begin{definition}\label{def:derived}
The bounded derived category $D^b(\mathcal{A})$ of $\mathcal{A}$ is obtained from $C^b(\mathcal{A})$ by inverting quasi-isomorphisms. 
The morphisms in $D^b(\mathcal{A})$ are formal compositions $f^{-1} \circ g$ where $f$ is a quasi-isomorphism. 
\end{definition}

\noindent
We denote with $\Stab(X)$ the set of all the Bridgeland stability conditions on $X$. In \cite{B2} Bridgeland proved that $\Stab(X)$ carries a natural geometric structure 
(see also \cite{BM1}, Theorem 2.2 and Proposition 2.3):

\begin{theorem} \label{th:stab-var}
$\Stab(X)$ is a finite dimensional complex variety. Moreover, given a fixed numerical invariant $v$, $\Stab(X)$ admits a locally finite dimensional walls-chambers decomposition. 
\end{theorem}

\noindent
Let $\mathcal{M}_{\sigma}(v)$ be the stack of $\sigma$-semistable objects in $D^b(\mathcal{A})$ of fixed numerical invariant $v$. Bayer and Macr\`i proved the following crucial result (see \cite{BM1}, Lemma 3.3):

\begin{lemma} \label{lem:posLemma} \emph{(Positivity Lemma)}
	Let $\sigma$ be a stability condition and $v$ a fixed numerical invariant. Then there exists a nef divisor $\mathcal{L}_{\sigma}$ on $\mathcal{M}_{\sigma}(v)$. Moreover, a curve $C$ on $\mathcal{M}_{\sigma}(v)$ is such that $\mathcal{L}_{\sigma}\cdot C = 0$ if and only if for two general closed points $c$ and $c'$ in $C$, the corresponding objects $\mathcal{E}_c$ and $\mathcal{E}_c' \in D^b(\mathcal{A})$ are $S$-equivalent.
\end{lemma}

\section{The ample cone of moduli spaces of sheaves on the plane}\label{amp}

\noindent
In this section we focus on the case of $\mbP^2$. In  \cite{B2} Bridgeland  proved that the space of stability conditions $\Stab(\mbP^2)$ contains as a subset the upper half plane $H \coloneq \{(s, t) | s \in \mbR, t > 0\}$. 
In this case, the numerical invariant $v$ is given by the Chern character $\xi$. 

\begin{theorem}\emph{\cite{ABCH}} \label{th:coarse}
	Let $\xi$ be a fixed Chern character. For each $(s, t) \in H$, there exists a coarse moduli space $M_{s,t}(\xi)$ parametrizing the $\sigma_{s,t}$-semistable objects with fixed Chern character $\xi$.
\end{theorem}

\noindent 
The geometry of $M_{s,t}(\xi)$ has been recently investigated by Li and Zhao in \cite{LZ2} (see also \cite{LZ1}).

\noindent
Fixing the Chern character $\xi$ implies, via Hirzebruch-Riemann-Roch, that we are also fixing the Euler characteristic and so the Hilbert polynomial. Therefore, we can consider the space $M(\xi)$
of $S$-equivalence classes of Gieseker semistable sheaves with Hilbert polynomial $P_\xi$.

\noindent
Bridgeland stability conditions have been proven to be useful to study the birational geometry of $M(\xi)$. In particular, the goal of this section is to compute the ample cone of $M(\xi)$, following the approach in 
\cite{CH}.

\begin{theorem}\emph{\cite{ABCH}} \label{thm:decoP2}
	$H \subseteq \Stab(\mbP^2)$ admits the following walls-chambers decomposition: there is a unique vertical wall and to the left of this wall there is a finite number of distinct nested semicircular walls. 
If we call the semicircular wall of maximal radius the Gieseker wall and the chamber outside this wall the Gieseker chamber, then for any $(s, t)$ outside the Gieseker wall $M_{s,t}(\xi)$ and $M(\xi)$ are isomorphic.
\end{theorem}

\noindent
Now we apply Bayer and Macr\`i Positivity Lemma to this situation.
\begin{lemma} \label{lem:posLemP2}
	For any $(s,t)$ in the Gieseker chamber, there exists a nef divisor $l_{s,t}$ on $M_{s,t}(\xi)$. Moreover, if we choose $(s, t)$ on the Gieseker wall, there exists a curve $C$ such that $l_{s, t} \cdot C = 0$, so $l_{s,t}$ is nef but not ample.
\end{lemma} 

\noindent
We go back to the problem of computing the ample cone of $M(\xi)$. As shown in Section 18 of \cite{LP}, the rank of the N\'eron-Severi of $M(\xi)$ can be only one or two. Moreover, as explained in \cite{CH}, Proposition 2.4, one of the extremal rays of the ample cone was already known in terms of the Donaldson-Uhlenbeck-Yau compactification. Thanks to Theorem \ref{thm:decoP2} and Lemma $\ref{lem:posLemP2}$, to determine the ample cone of $M(\xi)$ it is enough to compute the Gieseker wall.

\section{A first glimpse to formal definitions}\label{def}

\noindent
In this section we would like to give a gentle introduction to the notion of Bridgeland stability conditions, focusing again on the case of $\mbP^2$. 
The goal, more than precision, is to get the reader a little bit familiar with the notions that she will need to learn to enter in the technical core of the topic.
For more details, we refer to \cite{B} and \cite{H}.

\begin{definition} \label{def:TorsClass}
	A torsion pair in an abelian category $\mathcal{A}$ is a pair $(\mathcal{T}, \mathcal{F})$ of full additive subcategories of $\mathcal{A}$ such that 
	\begin{itemize}
		\item $\Hom_{\mathcal{A}}(\mathcal{T}, \mathcal{F}) = 0$ for all $T \in \mathcal{T}$ and $F \in \mathcal{F}$;
		\item Any $E \in \mathcal{A}$ fits into a short exact sequence 
		\begin{equation*}
		0 \to T \to E \to F \to 0,
		\end{equation*} 
		where $T \in \mathcal{T}$ is called \emph{torsion class} and $F \in \mathcal{F}$ is called \emph{torsion-free class}.
	\end{itemize}
\end{definition}

\begin{definition} \label{def:tStr}
	Let $D^b(\mathcal{A})$ be the bounded derived category of $\mathcal{A}$. A $t$-structure on $D^b(\mathcal{A})$ is a pair $(\mathcal{X}, \mathcal{Y})$ of full additive subcategories of $D^b(\mathcal{A})$ such that
	\begin{itemize}
		\item $\Hom_{D^b(\mathcal{A})}(X, Y) = 0$ for all $X \in \mathcal{X}$ and $Y \in \mathcal{Y}$;
		\item Any $D \in D^b(\mathcal{A})$ fits into a triangle 
		\begin{equation*}
		X \to D \to Y \to X[1] 
		\end{equation*} 
		where $X \in \mathcal{X}$ and $Y \in \mathcal{Y}$;
		\item $\mathcal{X}[1] \subseteq \mathcal{X}$. 
	\end{itemize}
The usual notation is $\mathcal{D}^{\leq 0} \coloneq \mathcal{X}$ and  $\mathcal{D}^{\geq 1} \coloneq \mathcal{Y}$. 
The heart of a $t$-structure is given by $\mathcal{H} \coloneq \mathcal{D}^{\leq 0} \cap \mathcal{D}^{\geq 0}$. The heart turns out to be an abelian category.
\end{definition}

\noindent
A $t$-structure is called bounded if every object in $D^b(\mathcal{A})$ is contained in $\mathcal{D}^{\leq n} \cap \mathcal{D}^{\geq -n}$ for $n$ sufficiently large.

\begin{remark} \label{rem:boundTstruStandard}
	Let $\mathcal{A}$ be an abelian category and let $D^b(\mathcal{A})$ be the bounded derived category of $\mathcal{A}$. Then, the standard bounded $t$-structure is given by 
	\begin{itemize}
		\item $\mathcal{D}^{\leq 0} = \{ X \in D^b(\mathcal{A}) | H^i(X) = 0 \textrm{ for } i > 0\}$;
		\item $\mathcal{D}^{\geq 0} = \{ Y \in D^b(\mathcal{A}) | H^i(Y) = 0 \textrm{ for } i < 0\}$.
		
	\end{itemize} 
\end{remark}

\noindent
Now we move on to the definition of Bridgeland stability conditions. We focus on the case of $\mbP^2$, so that $\mathcal{A} = \coh(\mbP^2)$ and $D^b(\mathcal{A}) = D^b(\coh(\mbP^2))$. 
In this part we mainly refer to \cite{ABCH}, Section 5.

\begin{definition} \label{def:stabCondition}
	A Bridgeland stability condition on $\mbP^2$ is a triple $(\mathcal{H}; r, d)$ such that 
	\begin{itemize}
	\item $\mathcal{H}$ is the heart of a bounded $t$-structure on $D^b(\coh(\mbP^2))$.
	\item $r$ and $d$ are linear maps 
	\begin{equation*}
	r, d : K(D^b(\coh(\mbP^2))) \to \mbR
	\end{equation*}
           defined on the $K$-group of the derived category, i.e. $r$ and $d$ are additive on triangles. The maps $r$ and $d$ satisfy 
	\begin{itemize}
		\item $r(E) \geq 0$ for each $E \in \mathcal{H}$;
		\item If $r(E) = 0$ for some non-zero object $E \in \mathcal{H}$, then $d(E) > 0$.
	\end{itemize}
	\item All the objects in $\mathcal{H}$ satisfy the property of Harder-Narasimhan.
	\end{itemize}
\end{definition}

\noindent
We now recall the concepts of slope, stability and Harder-Narasimhan filtration.

\begin{definition} \label{def:slope}
	The slope of a non-zero object $E \in \mathcal{H}$ with respect to $r, d$ is 
	\[
	\mu(E)=
	\begin{cases}
	d(E)/r(E) & \text{if $r(E) \neq 0$,} \\
	+ \infty & \text{otherwise.}
	\end{cases}
	\]
\end{definition}

\begin{definition} \label{def:stable}
	An object $E \in \mathcal{H}$ is called stable (resp. semistable) if for any proper non-zero subobject $F \subseteq E$
	\begin{equation*}
	\mu(F) < \mu(E) (\textrm{resp. } \leq)
	\end{equation*}
\end{definition}

\noindent
We recall that an object $E \in \mathcal{H}$ satisfies the property of Harder-Narasimhan if it admits a finite filtration
\begin{equation*}
0 \subseteq E_0 \subseteq E_1 \subseteq \dots \subseteq E_n = E
\end{equation*}
uniquely determined by the fact that $F_i \coloneq E_i/E_{i-1}$ is semistable and 
\begin{equation*}
\mu(F_1) > \mu(F_2) > \dots > \mu(F_n).
\end{equation*} 

\noindent
We now want to explicitly construct stability conditions on $\mbP^2$. We need, therefore, to exhibit a heart and the maps $r$ and $d$ such that they satisfy all the properties of Definition \ref{def:stabCondition}. 
The first choice could be the heart of the standard t-structure on $D^b(\coh(\mbP^2))$, see Remark \ref{rem:boundTstruStandard}, with the ``ordinary'' rank and degree, namely, for $E \in \mathcal{H}$:
\begin{itemize}
	\item $d(E) \coloneq c_1(E) \cdot L$;
	\item $r(E) \coloneq c_0 \cdot L^2$,
\end{itemize}
where $L$ is the hyperplane class on $\mathbb{P}^2$.
\noindent
Unluckily, these choices do not origin to a Bridgeland stability condition, because if we consider the skyscraper sheaf $\mbC_p$ we have
\begin{equation*}
r(\mbC_p) = 0 = d(\mbC_p).
\end{equation*}

\noindent
However, the associated slope function $\mu(E) = d(E)/r(E)$ satisfies the following weak property of Harder-Narasimhan: there exists a filtration 
\begin{equation*}
0 \subseteq E_0 \subseteq E_1 \subseteq \dots \subseteq E_n = E
\end{equation*}
\noindent
where $E_0$ is torsion-free and for $i > 0$, $F_i = E_i/E_{i-1}$ are torsion-free semistable sheaves with strictly decreasing slopes $\mu_i \coloneq \mu(F_i)$.

\noindent
Since the standard choice does not work, the idea is now to use the slopes $\mu_i$ to construct a torsion pair on $\coh(\mbP^2)$ inducing a bounded  $t$-structure on $D^b(\coh(\mbP^2))$.

\begin{definition}
	Let $s \in \mathbb{R}$ and let $(\mathcal{T}_s, \mathcal{F}_s)$ be full subcategories of $\coh(\mbP^2)$ such that 
	\begin{itemize}
		\item $T \in \mathcal{T}_s$ if $T$ is torsion or if every $\mu_i > s$ in the filtration of Harder-Narasimhan of $T$;
		\item $F \in \mathcal{F}_s$ if $F$ is torsion-free and if every $\mu_i \leq s$ in the filtration of Harder-Narasimhan of $F$.
	\end{itemize}
\end{definition}

\noindent
In this way, for any $s \in \mathbb{R}$, every pair $(\mathcal{T}_s, \mathcal{F}_s)$ is a torsion pair on $\coh(\mbP^2)$ that induces the following bounded $t$-structure on $D^b(\coh(\mbP^2)$:
\begin{itemize}
	\item $\mathcal{D}^{\geq 0} \coloneq \{\textrm{complexes } E | H^{-1}(E) \in \mathcal{F}_s \textrm{ and } H^i(E) = 0 \textrm{ for } i < -1\}$;
	\item $\mathcal{D}^{\leq 0} \coloneq \{ \textrm{complexes  } E | H^0(E) \in \mathcal{T}_s  \textrm{ and } H^i(E) = 0  \textrm{ for } i > 0\}$. 
\end{itemize}

\noindent
The heart of the $t$-structure defined by the torsion-pair is 
\begin{equation*}
H_s \coloneq  \{\textrm{complexes } E | H^{-1}(E) \in \mathcal{F}_s, H^0(E) \in \mathcal{T}_s  \textrm{ and } H^i(E) = 0 \textrm{ otherwise}\}. 
\end{equation*}

\noindent
In order to conclude, we just need to introduce the right notion of rank and degree.

\begin{theorem}\emph{\cite{B2}} \label{thm:stabcondP2}
	Let $s \in \mathbb{R}$ and $t > 0$. If we define the functions of rank and degree on $H_s$ in the following way:
	\begin{itemize}
		\item $r_t \coloneq t \cdot c_1(E(-s)) \cdot L$;
		\item $d_t \colon -(t^2/2) c_0(E(-2))\cdot L^2 + c_2(E(-s))$,
	\end{itemize}
then $(H_s; r_t, d_t)$ is a stability condition on $\mbP^2$ with slope function $\mu_{s,t} \coloneq d_t/r_t$.
\end{theorem}

\noindent
Indeed, it turns out that $r_t$ and $d_t$ satisfy all the properties of Definition \ref{def:stabCondition} (see \cite{ABCH}, remarks following Definition 5.10).

\end{document}